\documentclass[10pt]{article}

\usepackage{mathrsfs}
\usepackage{amsfonts}
\usepackage{amsmath}
\usepackage{amssymb}
\usepackage{amsthm}
\usepackage{latexsym}
\usepackage[active]{srcltx}
\usepackage{a4wide}

\usepackage{yhmath}

\def\R{{\mathbb R}}
\def\C{{\mathbb C}}
\def\N{{\mathbb N}}

\def\Z{{\mathbb Z}}

\def\Id{{\mathbb I}}

\def\P{{\mathcal P}}

\def\h{{\frak h}}

\def\OM{\Omega}

\def\om{\omega}

\def\ch{\chi}

\def\trans#1#2{{}^{#1}\mkern-3mu #2}

\def\Id{{\Bbb I}}

\def\iy{\infty}

\def\cds#1#2{#1,\cdots,#2}

\def \nn{\nonumber}

\def\res#1{_{\vert #1}}

\def\cd�?�??�?�

\def\cds{\cdots}

\def\lef({\left(}
\def\rig){\right)}

\def\noi{\noindent}

\let\e\varepsilon
\let\a\alpha

\newtheorem{theorem}{Theorem}
\newtheorem{prop}[theorem]{Proposition}
\newtheorem{cor}[theorem]{Corollary}
\newtheorem{lemma}[theorem]{Lemma}
\newtheorem{remarks}[theorem]{Remarks}
\newtheorem{remark}[theorem]{Remark}
\newtheorem{definition}[theorem]{Definition}
\newtheorem{proposition}[theorem]{Proposition}
\newtheorem{corollary}[theorem]{Corollary}
\newtheorem{notations}[theorem]{Notations}

\numberwithin{equation}{section}

\numberwithin{theorem}{section}

\newcommand{\blm}{\begin{lemma}}
\newcommand{\elm}{\end{lemma}}
\newcommand{\bpr}{\begin{proof}}
\newcommand{\epr}{\end{proof}}
\newcommand{\bcor}{\begin{corollary}}
\newcommand{\ecor}{\end{corollary}}

\begin{document}

\title
{Harmonic analysis of weighted $L^p$-algebras}

\author{Yu. N. Kuznetsova \and C. Molitor-Braun
\thanks{Supported by the
research grant F1R-MTH-PUL-10NCHA of the University of
Luxembourg}}

\date{}

\maketitle

\begin{abstract} Let $G$ be a locally compact, compactly generated
group of polynomial growth and let $\omega$ be a weight on $G$.
Under proper assumptions on the weight $\omega$, the Banach space
$L^p(G,\omega)$ is a Banach $*$-algebra. In this paper we give
examples of such weighted $L^p$-algebras and we study some of
their harmonic analysis properties, such as symmetry, existence of
functional calculus, regularity, weak Wiener property, Wiener
property, existence of minimal ideals of a given hull.
 \footnote{keywords: Weighted $L^p$-algebras, symmetry, functional
 calculus, regularity, weak Wiener property, Wiener property,
 minimal ideals of a given hull
-- 2010 Mathematics Subject Classification:
43A15, 
22D15, 
22D20  
}
\end{abstract}

\section{Introduction}

Weights and weighted function spaces play an important role in
mathematics. In essence, a weight makes it possible to study the
behaviour of functions around a certain point, ignoring their
oscillations at infinity, or on the contrary, to amplify the
asymptotic behaviour of a function. More precisely, introducing a
weight means modelling in a quantitative manner the decay of the
functions to be studied. This has numerous applications in
numerical mathematics and is quite often used for concrete
applications (signal theory, Gabor analysis, sampling theory,...),
see for instance (\cite{Groe.Lei.}, \cite{Da.Fo.Groe.},
\cite{Groe.Lei.1}, \cite{Fe.Groe.Lei.}).

\noi On the other hand, weights appear naturally in analysis: in
inequalities relating the norm of a function to the norm of its
derivatives, in extension theorems, etc.; see, e.g., a survey of
L.~D.~Kudryavtsev and S.~M.~Nikol'sky \cite{kud-nik}. One of the
areas where weighed spaces are applied most intensively is the
theory of boundary value problems for partial differential
equations (see the surveys \cite{kud-nik}, \cite{glushko}). By
the way, using the Laplace transform also means working in a weighted
function space.

\noi In representation theory, which interests us most, weights
occur for instance in the following way. If $G$ denotes a locally
compact group and $(T,V)$ is a continuous representation of $G$,
then the maps $\omega:x \mapsto \Vert T(x)\Vert_{op}$ and
$\omega:x\mapsto \max(\Vert T(x)\Vert_{op}, \Vert
T(x^{-1}\Vert_{op})$ are weights, the last one being symmetric
($\omega(x^{-1})=\omega(x),~\forall x$). For any one of these
weights, the map
$$f \mapsto T(f):=\int_G f(x)T(x)dx$$ is a representation of the
weighted function algebra $$L^1(G,\omega):= \{f:G\rightarrow
\C~\vert~f \text{ measurable and } \int_G \vert f(x)\vert
\omega(x)dx < +\infty\}.$$ Subexponential weights like the ones
introduced in \ref{subexp-w} below, appear in the context of
nilpotent Lie groups. In fact, let $G$ be a connected, nilpotent
Lie group. Let $G_1$ be the derived group of $G$, i. e. the closed
subgroup generated by the elements of the form $[x,y]=
x^{-1}y^{-1}xy,~x,y\in G$. Let $U$ be a generating neighbourhood
of the identity $e$ in $G$ and $V=U\cap G_1$ the corresponding
neighbourhood of $e$ in $G_1$. Let
$$\vert x \vert_U:= \inf\{n \in \N~\vert~x\in U^n\}~\text{ for
}x\neq e,~\vert e \vert_U=1,$$ similarly for $\vert x\vert_V$.
Then it is shown in \cite{Al.}), that for any weight $\omega$ on
$G$ which is submultiplicative, i. e. such that $\omega(xy) \leq
\omega(x)\omega(y)$ for all $x,y$,
$$\omega\vert_{G_1}(x) \leq e^{C \vert x
\vert_V^{\frac{1}{2}}},~\forall x\in G_1,$$ for some constant $C$.
By the way, on any compactly generated locally compact group, with
generating neighbourhood $U$, every submultiplicative weight
$\omega$ is exponentially bounded, i. e. satisfies a relation of
the form $\omega(x) \leq e^{K\vert x \vert_U}$ for $K=\ln
\sup_{x\in U}\omega(x)$.

\noi If the weight $\omega$ is submultiplicative, then the
weighted function space $L^1(G,\omega)$ is a Banach algebra for
convolution, and even a Banach $*$-algebra if the weight is
symmetric. The advantage of Banach $*$-algebras over just Banach
spaces is clear. They have a much richer structure which may be
studied via representation theory and harmonic analysis
techniques. In this way, interesting problems arise. Let us just
mention the question of their ideal theory, problems of
generalized spectral synthesis, of symmetry of the algebra, of
invertibility, factorization problems. All these questions make
sense for the weighted algebra $L^1(G,\omega)$. Moreover, it is
the harmonic analysis properties for algebras, that make the
weighted algebras $L^1(G,\omega)$ interesting for some of the
concrete applications mentioned in the beginning (see for instance
\cite{Groe.Lei.1}).

\noi On the other hand, the importance of $L^p$-spaces of the form
$L^p(G)$ or $L^p(G,\omega)$, $1 < p < +\infty$, is well known in
functional analysis. It would be attractive to extend the theory
of convolution algebras to the $L^p$-case, because $L^p$ spaces
are reflexive --- not a common property among Banach algebras.
Unfortunately, if $G$ is not compact and if $p\neq 1$, $L^p(G)$ is
not an algebra for convolution. Nevertheless, for appropriate
groups $G$ and weights $\omega$, the weighted $L^p$-spaces
$$L^p(G,\omega):=\{ f:G\rightarrow\C~\vert~f\text{ measurable and
} \Vert f \Vert_{p,\omega}:=\big ( \int_G \vert
f(x)\vert^p\omega(x)^pdx \big)^{\frac{1}{p}} < +\infty \}$$ may be
algebras. A sufficient condition on the group $G$ for the
existence of weighted $L^p$-algebras, is that the group is
$\sigma$-compact. In that case, there are even a lot of such
weighted $L^p$-algebras. In \ref{product-w} we show that any
positive symmetric submultiplicative function multiplied by any
$L^p$-algebra weight produces again an $L^p$-algebra weight. This
makes it possible to construct $L^p$-algebra weights with all
kinds of different growth behaviors.

\noi In the context of weighted $L^p$-algebras, let us mention the
works of Wermer \cite{wermer}, Nikol'ski \cite{nik-mian},
Feichtinger \cite{Fei}, and recently (\cite{Ku.1}-\cite{Ku.4}).
Most of these papers concentrate mainly on the question whether
the corresponding $L^p$-spaces are algebras. The only well-studied
case is $L^p(\Z,\omega)$, see, e.g., a long paper of El-Fallah,
Nikol'ski and Zarrabi \cite{el-fallah}. This is mainly for the
reason that in the problems of weighted approximation by
polynomials, as initiated by S.~N.~Bernstein \cite{bernstein},
$L^p(\Z, \omega)$ algebras play a distinguished role
\cite{nik-mian}. But this is not the only possible application of
weighted $L^p$-algebras. Similarly as for $L^1$-algebras, the
weights can be used in numerical mathematics to model the decay of
the functions to be used and allow numerical computations. On the
other hand, weighted $L^p$-algebras may turn out to be important
examples for people working in Banach algebra theory or operator
theory. Such algebras have already been used successfully in the
interpolation theory and in questions of factorization
\cite{Bl.Ka.Ra.}.

\noi As for possible harmonic analysis properties of an arbitrary
Banach $*$-algebra ${\mathcal A}$, several questions have
particularly caught the interest of mathematicians:

\noi Is the algebra symmetric, i. e. does every self-adjoint
element have a real spectrum?

\noi Is the algebra ${\mathcal A}$ regular, i. e. do the elements
of the algebra separate points from closed sets in
$\text{Prim}_*{\mathcal A}$, the space of kernels of topologically
irreducible unitary representations?

\noi Does the algebra have the weak Wiener property, i. e. is
every proper, closed, two-sided ideal annihilated by an
algebraically irreducible representation?

\noi Does the algebra have the Wiener property, i. e. does the
previous property hold for topologically irreducible unitary
representations?

\noi Do there exist minimal ideals of a given hull?

\noi The list of authors who studied group algebras $L^1(G)$ and
their properties is long and is outside the scope of this paper.
For weighted group algebras $L^1(G,\omega)$, let us mention among
others the following: In the abelian case, the systematic study of
such properties for weighted group algebras $L^1(G,\omega)$ goes
back to Beurling (\cite{Beu.1}, \cite{Beu.2}), Domar \cite{Do.}
and Vretblad \cite{Vr.} among others. In the non-abelian case, one
may refer to Hulanicki (\cite{Hu.1}, \cite{Hu.2}), Pytlik
(\cite{Py.1}, \cite{Py.2}), as well as to more recent studies
(\cite{Dz.Lu.Mo.}), \cite{Fe.Groe.Lei.Lu.Mo.} and
\cite{Fe.Groe.Lei.}). In \cite{Fe.Groe.Lei.} for instance, the
question of the symmetry for weighted group algebras
$L^1(G,\omega)$ is completely solved for compactly generated
groups with polynomial growth. Let us mention in this respect,
that these abstract problems may be quite important for concrete
applications. Hence Gr\" ochenig and Leinert \cite{Groe.Lei.}
point out that the theory of symmetric group algebras is an
important tool to solve problems about Gabor frames, motivated by
signal theory.

\noi A systematic study of harmonic analysis properties of
weighted $L^p$-algebras $L^p(G,\omega)$ should also be of
importance, as well for applications and for more abstract
mathematical problems.  In \cite{Ku.3} some harmonic analysis
properties like the regularity are studied in the case of abelian
groups. But not much seems to have been done up to now in the
non-abelian case. Hence the main purpose of this paper will be to
study harmonic analysis properties in the context of non-abelian
weighted $L^p$-algebras.

\noi In the present paper, we work on general compactly generated
groups with polynomial growth. The weight $\omega$ is supposed to
satisfy $\omega^{-q}*\omega^{-q} \leq C \omega^{-q}$ for some
constant $C>0$, where $\frac{1}{p}+ \frac{1}{q}=1$. This ensures
$L^p(G,\omega)$ to be an algebra (\cite{wermer}, \cite{Ku.1}). We
also assume the weight to be submultiplicative. We start by giving
examples of $L^p(G,\omega)$ $*$-algebras, as well for polynomially
growing weights as for sub-exponentially growing weights. We then
address the questions raised previously: We prove the symmetry of
the $L^p$-algebra $L^p(G,\omega)$, if either $G$ is abelian and
$\omega$ satisfies the same condition as for the case of
$L^1(G,\omega)$ (condition (S), \cite{Fe.Groe.Lei.Lu.Mo.},
\cite{Fe.Groe.Lei.}) or if $\omega$ is polynomial in the sense of
Pytlik \cite{Py.2} (and $G$ not necessarily abelian). The same
hypothesis as in \cite{Dz.Lu.Mo.}, i. e. the non-abelian
Beurling-Domar condition (BDna) allows us to construct a
functional calculus on a total subset of the algebra
$L^p(G,\omega)$ and to show regularity, as well as the weak Wiener
property. If the $L^p$-algebra is moreover symmetric, we also get
the Wiener property and the existence of minimal ideals of a given
hull. Let us recall that the (BDna) condition is defined as
follows in \cite{Dz.Lu.Mo.}: Let $G=\bigcup_n U^n$, where $U$ is a
relatively compact, generating neighbourhood of $e$ in $G$. We
define $s(n):=\sup_{x\in U^n} \omega(x)$. Then the weight $\omega$
satisfies (BDna) if and only if $$\sum_{n\in \N, n\geq e^e} \frac{
\big ( \ln(\ln n) \big ) \ln \big ( s(n) \big)}{1+n^2} <
+\infty.$$ This condition is independent of the choice of the
generating neighbourhood $U$ and is only slightly stronger than
the conditions used by Domar and Beurling in the abelian case (see
\cite{Dz.Lu.Mo.} for additional comments). These results relay on
the corresponding results for $L^1$-algebras $L^1(G,\omega)$. The
question of whether, more generally, condition (S) or the
GRS-condition as defined in \cite{Fe.Groe.Lei.Lu.Mo.} and
\cite{Fe.Groe.Lei.} imply the symmetry of the algebra
$L^p(G,\omega)$, as they do for $L^1(G,\omega)$, is still an open
problem. Finally, let us point out that although these algebras
$L^p(G,\omega)$ have certain nice harmonic analysis properties,
they are not amenable if $p>1$ and $G$ non-discrete, as they don't
have bounded approximate identities \cite{Ku.2}.

\noi Let us also mention, that we assume our weights to be
submultiplicative, in order for $L^1(G,\omega)$ to be a
$*$-algebra. Therefore we can rely on the known
$L^1(G,\omega)$-results. But there are examples of weights which
are not submultiplicative and which produce nevertheless Banach
$*$-algebras $L^p(G,\omega)$ (see for instance \cite{Ku.2}).
Studying these weights and the properties of the corresponding
$L^p$-algebras is still a challenge.

\subsection{Assumptions on groups and weights}

We suppose in this paper that $G$ is a compactly generated locally compact group. This group $G$ is a group of \emph{polynomial growth} if there is a relatively compact generating neighbourhood of the identity $U$ such that $|U^n|\le C n^Q$ for some constants $C,Q$. It is known that $Q$ does not depend on the choice of $U$. The class of such groups will be denoted by [PG], and for the rest of the paper we assume that $G$ is [PG].

\noi If $U$ is a relatively compact generating neighbourhood of the identity, we define
$$|x|_U := \inf\{ n~\vert~x\in U^n\}.$$
When the choice of $U$ is not important, we write simply $|x|$. In the case $G=\R$, $|x|$ may also denote the absolute value of $x$, and the results of this paper remain correct.

\noi Let $$\omega: G \rightarrow [1, +\infty[ $$ be a measurable function (\emph{weight}) such that
\begin{eqnarray}\label{}
\nn \omega(xy) &\leq& \omega(x) \omega(y), \quad \forall x,y \in G\\
\nn \omega(x)&=&\omega(x^{-1}), \quad \forall x\in G\\
\nn \omega^{-q} * \omega^{-q} &\leq& \omega^{-q},
\end{eqnarray}
where $\frac{1}{p}+\frac{1}{q}=1$, $p > 1$. These conditions are
sufficient for $L^p(G,\omega)$ to be a $*$-Banach algebra
(\cite{wermer}, \cite{Ku.1}). We will say that $(G,\omega)$
satisfies (LPAlg) ($L^p$-algebra) if these conditions are
satisfied.

\noi It is often easier to check that $\omega$ satisfies
conditions (LPAlg) with some constants $C_1,C_2$: $\omega(xy) \leq
C_1\omega(x) \omega(y)$ and $\omega^{-q} * \omega^{-q} \leq
C_2\omega^{-q}$. But a renormalizing $\omega_1=C\omega$ with
$C=\max(C_1,C_2^{1/q})$ gives an equivalent weight satisfying
(LPAlg).

\noi Since every group in [PG] is amenable, it follows from
(LPAlg), by \cite[Theorem~3.2]{Ku.4}, that $\omega^{-q}\in
L^1(G)$. We may assume, without loss of generality, that the
weight $\omega$ is continuous (\cite{Fei}). This will be assumed
for the rest of the paper, except for some examples which depend
on the discontinuous function $|x|=\vert x\vert_U$.

\subsection{Examples of weights}

\subsubsection{Polynomial weights}\label{poly-w}

On every group of polynomial growth, the weight $\omega(x)=(1+|x|)^D$ satisfies (LPAlg) for $D$ sufficiently large (\cite{Fei}), so $L^p(G,\omega)$ is an algebra.

\subsubsection{Products of weights}\label{product-w}

Let $u$ be a positive submultiplicative function on $G$, $u(x)= u(x^{-1})$, and let $w_1$ be a weight satisfying (LpAlg). Then $w(x) = u(x) w_1(x)$ also satisfies (LpAlg). In particular, any such submultiplicative function $u$, multiplied by $(1+|x|)^D$ for $D$ sufficiently large, is a (LpAlg)-weight. To prove this, we need to check only the last condition:

$$
(w^{-q}*w^{-q}) (x) = \int_G w^{-q}(y) w^{-q}(y^{-1}x) dy
 = \int_G  u^{-q}(y) u^{-q}(y^{-1}x) w_1^{-q}(y) w_1^{-q}(y^{-1}x) dy.
$$
From $u(x)\le u(y)u(y^{-1}x)$ we have $u(x)^{-q}\ge u(y)^{-q}u(y^{-1}x)^{-q}$, so the integral above is bounded by
$$
(w^{-q}*w^{-q}) (x)
 \le  u^{-q}(x) \int_G  w_1^{-q}(y) w_1^{-q}(y^{-1}x) dy
 \le  u^{-q}(x) w_1^{-q}(x) = w^{-q}(x).
$$

\subsubsection{Non-polynomial weight}\label{subexp-w}

By the reasoning of section \ref{product-w},
$\omega(x)=e^{|x|^\gamma}(1+|x|)^D$ with $0<\gamma\le 1$ is an
$L^p$-algebra weight on any group in [PG] for all $D$ sufficiently
large.

\noi Moreover, it can be shown that the weight $\omega(x)=e^{|x|^\gamma}$ itself with $0<\gamma<1$ satisfies (LPAlg) for all $p>1$. For $G=\R$ this example is contained in the very first paper of Wermer \cite{wermer} on $L^p$-algebras.

\noi Let $G=\cup U^n$, where $U=U^{-1}$ and $|U^n|\le Cn^Q$. Denote $U_n=U^n\setminus U^{n-1}$ and assume that $G$ is non-compact. We define
$$
\omega_n\equiv \omega|_{U_n} = e^{n^\gamma},
$$
$0<\gamma<1$, and show that $L^p(G,\omega)$ is an algebra for every $p>1$ (the case $p=1$ is known). For this, we check the sufficient condition $\omega^{-q}*\omega^{-q}\le C' \omega^{-q}$.

\noi Denote $u=\omega^{-q}$, then
$$
u_n\equiv u|_{U_n} = e^{-qn^\gamma}.
$$
Take $x\in U_m$ and estimate $\dfrac{u*u}{u}(x)$.
$$
\dfrac{u*u}{u}(x) = {1\over u_m}\int_G u(y)u(y^{-1}x)dy = \sum_n {u_n\over u_m} \int_{U_n} u(y^{-1}x)dy.
$$
If $y\in U_n=U_n^{-1}$, and $y^{-1}x\in U_k$, then $\max( n-m, m-n) \le k\le n+m$.
Denote $U^x_{nk}=U_k\cap (U_nx)$. Then $U_nx=\cup_k U^x_{nk}$, $U_k=\cup_n U^x_{nk}$. In particular, $|U^x_{nk}|\le \min( |U_n|,|U_k|)$. We can rewrite the sum as
$$
\dfrac{u*u}{u}(x) = \sum_{n,k} {u_n\over u_m} \,u_k |U^x_{nk}|.
$$
Note that $u_n$ decreases and $C_0=\int_G u<\infty$. Split now the sum into four parts:

1) $n\ge m$. Then $u_n\le u_m$, and
$$
({\rm sum1}) \le \sum_{n,k} {u_m\over u_m} \,u_k |U^x_{nk}|\le \sum_k u_k|U_k| = \int_G u=C_0.
$$

2) Similarly, if $n<m$ but $k\ge m$ then $u_k\le u_m$ and
$$
({\rm sum2}) \le \sum_{n,k} {u_m\over u_m} \,u_n |U^x_{nk}|\le \sum_n u_n|U_n| = \int_G u=C_0.
$$

3) $m/2< n< m$, $m/2< k<m$. Then $\max(u_n,u_k)\le u_{[m/2]+1}\le \exp(-qm^\gamma/2^\gamma)$;
\begin{align*}
({\rm sum3}) &\le \sum_{n=[m/2]+1}^m\sum_{k=[m/2]+1}^m e^{q(m^\gamma-2m^\gamma2^{-\gamma})} |U^x_{n,k}|
 \\ {}&\le \sum_{n=1}^m e^{qm^\gamma(1-2^{1-\gamma})} |U_n|
 \le e^{qm^\gamma(1-2^{1-\gamma})} Cm^Q\cdot m.
\end{align*}
Since $2^{1-\gamma}>1$, the coefficient in the exponent is negative, so this expression tends to zero as $m\to\infty$. Thus, this is bounded by a constant $C_1$.

4) The only complicated case is $\,n\le m/2$, $m-n\le k<m$, and the symmetric case $k\le m/2$, $m-k\le n<m$ which reduces to the first one by exchanging $k$ and $n$. Here $u_k\le u_{m-n}$, so
\begin{align*}
({\rm sum 4}) &\le \sum_{n=0}^{[m/2]}{u_nu_{m-n}\over u_m}\sum_{k=m-n}^m  |U^x_{nk}|
\le \sum_{n=0}^{[m/2]} {u_nu_{m-n}\over u_m}|U_n| \\
&\le C\sum_{n=0}^{[m/2]}  e^{q(m^\gamma-n^\gamma -(m-n)^\gamma)} n^Q.
\end{align*}
\noi Denote
$$
f(x) = \int_0^{x/2}t^Q  e^{q(x^\gamma-t^\gamma -(x-t)^\gamma)} dt;
$$
clearly our sum is bounded for all $m$ if and only if $f$ is bounded on $\R_+$.

\noi By changing the variable to $s=t/x$ we have:
$$
f(x) = \int_0^{1/2} x^Q s^Q e^{qx^\gamma(1-s^\gamma -(1-s)^\gamma)} xds = x^{Q+1} \int_0^{1/2} t^Q e^{qx^\gamma(1-t^\gamma -(1-t)^\gamma)} dt.
$$
There is a classical theorem \cite[\S2.4]{erd} for integrals of the type
$$
F(x) = \int_\a^\beta g(t) e^{xh(t)}dt.
$$
Suppose that:
\begin{itemize}\renewcommand{\labelitemi}{}
\item $h$ is real-valued and continuous at $t=\a$;
\item $h'$ exists and is continuous for $\a<t\le\beta$;
\item $h'<0$ for $\a<t<\a+\eta$ with some $\eta>0$;
\item $h(t)\le h(\a)-\e$ for some $\e>0$ and all $t\in[\a+\eta,\beta]$;
\item $h'(t)\sim -a(t-\a)^{\nu-1}$ as $t\to \a$, where $\nu>0$;
\item $g(t)\sim b(t-\a)^{\lambda-1}$ as $t\to\a$, where $\lambda>0$.
\end{itemize}
Then
$$
f(x)\sim {b\over\nu} \Gamma\Big({\lambda\over\nu}\Big) \left({\nu\over ax}\right)^{\lambda/\nu} e^{xh(\a)}
$$
as $x\to\infty$.

\noi We can apply this theorem to
$$
F(x) = \int_0^{1/2} t^Q e^{x(1-t^\gamma -(1-t)^\gamma)} dt.
$$
In this case $g(t)=t^Q$, $h(t)=1-t^\gamma-(1-t)^\gamma$, $\a=0$, $\beta=1/2$. The derivative $h'(t)=-\gamma\big(t^{\gamma-1}-(1-t)^{\gamma-1}\big)$ is negative on $(0,1/2)$ since $\gamma-1<0$ (so that $t^{\gamma-1}$ decreases) and $t<1-t$. It follows that $h$ decreases, so all the conditions hold. We have $b=1$, $\lambda=Q+1$, $a=\nu=\gamma$. Thus,
$$
F(x) \sim {1\over\gamma} \Gamma\Big({Q+1\over\gamma}\Big) \left({\gamma\over \gamma x}\right)^{(Q+1)/\gamma} e^{x\cdot 0}
 \equiv C_2 x^{-(Q+1)/\gamma}.
$$
If we return to $f$, then we get
$$
f(x) = x^{Q+1} F(qx^\gamma) \sim C_2 x^{Q+1} (qx^\gamma)^{-(Q+1)/\gamma} = C_2 q^{-(Q+1)/\gamma} \equiv C_3.
$$
It follows that there is a constant $C_4$ such that $f(x)\le C_4$ for all $x>0$.

\noi Now, collecting all together, we have
$$
{u*u\over u}(x) \le 2C_0+C_1+C_4\equiv C',
$$
what completes the proof.

\subsubsection{Fast-growing weights}

If the weight is submultiplicative, as we always assume, then it can grow at most exponentially. But $L^p(G,\omega)$ is in general not an algebra with the weight $\omega(x)=e^{|x|}$. We will show this for $G=\R$. Take nonnegative $f,g\in L^p(\R,e^{|x|})$, then $F=e^{|x|}f,G=e^{|x|}g$ are in $L^p(\R)$;
$$
(f*g)(s) = \int_{-\infty}^\infty f(t)g(s-t)dt \ge \int_0^s F(t) e^{-t} G(s-t) e^{-s+t}dt = e^{-s} \int_0^s F(t) G(s-t) dt.
$$
Let $F_+=F\cdot I_{[0,+\infty)}$, $G_+=G\cdot I_{[0,+\infty)}$,
where $I_{[0,+\infty)}$ is the characteristic function of the
interval $[0,+\infty)$. If we assume that $f*g\in L^p(\R,\omega)$,
then from the formula above $F_+*G_+\in L^p(\R)$. Since for every
$F_+,G_+\in L^p([0,+\infty))$ we have $|e^{-t}F_+|, |e^{-t}G_+|\in
L^p(\R,\omega)$, it follows that $L^p([0,+\infty))$ is a
convolution algebra if $L^p(\R,\omega)$ is so. But it is
well-known that this is not true, if $p\neq 1$.

\noi Nevertheless, by \ref{poly-w} and \ref{product-w},
$L^p(\R,\omega_1)$ is an algebra for $\omega_1(x)=(1+\vert
x\vert)^D e^{\vert x\vert}$.

\noi There are, moreover, super-exponential weights with which
$L^p(G,\omega)$ is an algebra. One example is $\omega(x)=e^{x^2}$
on the real line, found by El Kinani \cite{elkin}. But this weight
is not submultiplicative.

\section{First properties of weighted algebras}

\subsection{Known inequalities}

The conditions (LPAlg) guarantee that $L^p(G,\omega)$ is an
algebra, that $L^p(G,\omega)$ is translation-invariant and
$\|\trans xf\|_{p,\omega} \le \omega(x)\|f\|_{p,\omega}$, where
$\trans xf(t)=f(x^{-1}t)$, as
$$
\|\trans xf\|_{p,\omega}^p = \int |f(x^{-1}t)|^p \omega(t)^pdt = \int |f(y)|^p \omega(xy)^pdy \le
 \omega(x)^p \|f\|_{p,\omega}^p.
$$

\noi Also under (LPAlg) the following is known:
\begin{itemize}
\item $L^1(G,\omega)\subset L^1(G)$, and $\|f\|_1\le \|f\|_{1,\omega}$ for all $f\in L^1(G,\omega)$
\item $L^p(G,\omega)\subset L^1(G)$, and $\|f\|_1\le C\|f\|_{p,\omega}$ for all $f\in L^p(G,\omega)$, with $C=\big(\int \omega^{-q}(x) dx\big)^{1\!/\!q}$
\item $L^1(G,\omega)*L^p(G,\omega)\subset L^p(G,\omega)$, and
\begin{equation}\label{Lp-L1-module}
\|f*g\|_{p,\omega}\le \|f\|_{1,\omega}\|g\|_{p,\omega}
\qquad\text{ for all } f\in L^1(G,\omega), ~g\in L^p(G,\omega)
\end{equation}
\end{itemize}

\noi If $\omega(x)=\omega(x^{-1})$, the usual involution $f^*(x)=\overline{f(x^{-1})}$ is an isometry on $L^p(G,\omega)$: $\|f^*\|_{p,\omega} = \|f\|_{p,\omega}$ (recall that every group of polynomial growth is unimodular).

\subsection{Approximate units}\label{sec-approx}

Under the assumption that $\omega$ is continuous, the proof of
(\cite{Ku.2}) of the fact that the measurable, bounded functions
of compact support are dense in $L^p(G,\omega)$, shows that the
same is also true for the set ${\mathcal C}_c(G)$ of continuous
functions with compact support. By (\cite{Ku.2}), there exists a
net $(f_s)_s$ of measurable, bounded functions with compact
support which form a bounded approximate identity in
$L^1(G,\omega)$ and an (unbounded) approximate identity in
$L^p(G,\omega)$. In fact, in the same way it can be proved that if $V_s$ runs through a basis of compact, symmetric neighbourhoods of the identity $e$ in $G$, then every family $f_s$ such that $0\le f_s\le1$, $\Vert f_s\Vert_1=1$, $\text{supp}f_s \subset V_s$, is an approximate identity with properties as above. It is easy to see that these functions $f_s$ may be chosen to be continuous and self-adjoint, $f_s=f_s^*$. Moreover, it will be convenient to have $V_s\subset K$, where $K$ is a fixed compact set.

\subsection{Representations of the weighted algebras}

For every non-degenerate $*$-representation $T$ of $L^p(G,\omega)$ in a Hilbert space $\cal H$, there is a unitary continuous representation $V$ of $G$ such that
\begin{equation}\label{int-TV}
T(f)=\int_G V(x)f(x)dx
\end{equation}
for all $f$. This is proved exactly like in \cite[Theorem~22.7]{HR}, though the assumption (ii) of this theorem does not hold. Let us make the following remarks on the proof.

\noi (1) If $T$ is cyclic with the cyclic vector $\xi$, one defines $V$ as follows: on the dense subspace of vectors of the type $T(f)\xi$, $f\in L^p(G,\omega)$, put $V(x)\big(T(f) \xi\big)=T(\trans xf)\xi$; it may be easily shown that this is an isometry, so $V(x)$ extends to a unitary operator on $\cal H$. It is also straightforward that $V$ is a representation.

\noi (2) To get the equality \eqref{int-TV}, we need to prove that
every coefficient $x\mapsto\langle V(x) T(f)\xi,\eta\rangle$,
where $f\in L^p(G,\omega)$, and $\xi,\eta\in \cal H$, is
measurable (this is \cite[22.3i]{HR}). But we have even more:
coefficients of $V$ are continuous since by \cite{Ku.2} the
mapping $x\mapsto \trans xf$ from $G$ to $L^p(G,\omega)$ is
continuous for every $f$. From this, by \cite[22.3]{HR} we get
some representation $\tilde T$ of $L^p(G,\omega)$ defined by
\begin{equation}\label{TV}
\tilde T(f) = \int_G V(x) f(x)dx
\end{equation}

\noi (3) To prove that $\tilde T=T$, take a vector $\eta\in \cal H$ and the cyclic vector $\xi \in \cal H$. For the linear functional $H(f)= \langle T(f) \xi,\eta\rangle$ on $L^p(G,\omega)$ there is a function $h\in L^q(G,\omega)$ such that $H(f) = \int f(x) h(x) dx$. Then
for any $f,g\in L^p(G,\omega)$ we have, with all integrals absolutely converging,
\begin{align*}
\langle \tilde T(f) T(g) \xi,\eta\rangle &= \int_G \langle V(x) f(x) T(g)\xi,\eta \rangle dx
 = \int_G \langle V(x) T(g)\xi,\eta \rangle f(x) dx
 \\&=
 \int_G \langle T(\trans xg)\xi,\eta \rangle f(x) dx
 = \int_G \left( \int_G \trans xg(y) h(y) dy \right) f(x) dx
\\& = \int_G \left( \int_G g(x^{-1}y) h(y) dy \right) f(x) dx
  = \int_G h(y) (f*g)(y) dy
\\&= \langle T(f*g)\xi,\eta\rangle
  = \langle T(f) T(g) \xi,\eta\rangle.
\end{align*}
It follows that $\tilde T=T$ on the dense subspace of vectors of the type $T(f)\xi$, $f\in L^p(G,\omega)$, and as a consequence on the whole $\cal H$.

\noi (4) In general, $T$ can be expanded into a direct sum of cyclic representations $T_\alpha$ \cite[21.13]{HR}. Every $T_\alpha$ is given by the formula \eqref{TV} with some $V_\alpha$; then $T$ is equal to the same integral \eqref{TV} with $V=\oplus V_\alpha$.

\noi Further, by \cite[22.6]{HR}, $T$ and $V$ are irreducible or not simultaneously.

\noi In particular, this gives us the identification
$\widehat{L^{\!p}(G,\omega)} = \widehat G$.

\section{Symmetry}

\subsection{ }
The notion of symmetry plays an important role in the theory of
Banach $*$-algebras. It may be defined as follows:

\noi Let $\mathcal A$ be a Banach $*$-algebra and let $a\in
{\mathcal A}$. We will denote the spectrum of $a$ in ${\mathcal
A}$ by $\sigma_{\mathcal A}(a)$ and the spectral radius of $a$ in
$\mathcal A$ by $r_{\mathcal A}(a)$. Then the algebra ${\mathcal
A}$ is said to be {\it symmetric} if $\sigma_{\mathcal A}(a^*a)
\subset [0, +\infty[$ for all $a\in {\mathcal A}$, or,
equivalently, if $\sigma_{\mathcal A} (a) \subset \R$ for all
$a=a^* \in {\mathcal A}$.

\noi For abelian Banach $*$-algebras the symmetry is equivalent to
the fact that all the characters of the algebra are unitary.

\noi Let $(G,\omega)$ satisfy (LPAlg). Let us recall the following definitions for
the weight $\omega$ (see for instance \cite{Fe.Groe.Lei.Lu.Mo.}
and \cite{Fe.Groe.Lei.}).

\begin{definition}
a) The weight $\omega$ on $G$ is said to satisfy the GRS-condition
(or GNR-condition), if
\begin{equation}\tag{GRS}
\lim_{n\rightarrow+\infty} \omega(x^n)^{\frac{1}{n}}=1, \quad
\forall x\in G.
\end{equation}

\noi b) The weight $\omega$ is said to satisfy condition (S), if,
for every generating, relatively compact neighbourhood $U$ of $G$,
\begin{equation}\tag{S}
\lim_{n\rightarrow+\infty} \sup_{x\in U^n}
\omega(x)^{\frac{1}{n}}=1.
\end{equation}
\end{definition}

\noi These conditions are linked to the symmetry of weighted group
algebras. Among others, the following results are known:

\medskip \noi For $G=\Z$, $l^1(\Z,\omega)$ is symmetric if and only if
$\lim_{n\rightarrow+\infty} \omega(n)^{\frac{1}{n}}=1$
(\cite{Nai.}).

\noi If $G\in [PG]$, then $L^1(G)$ is symmetric (\cite{Lo.}).

\noi If $G\in [PG]$ and $\omega$ satisfies condition (S), then
$L^1(G,\omega)$ is symmetric (\cite{Fe.Groe.Lei.Lu.Mo.}).

\medskip \noi The final version of results of this type is due to
Fendler, Gr\" ochenig and Leinert (\cite{Fe.Groe.Lei.}). They
prove:

\begin{theorem}\label{symmetry-l1}
Let $G\in [PG]$ and let $\omega$ be a weight on $G$. Then the
following are equivalent:

\noi (i) $\omega$ satisfies the GRS-condition.

\noi (ii) $\omega$ satisfies condition (S).

\noi (iii) $L^1(G,\omega)$ is symmetric.

\noi (iv) $\sigma_{L^1(G,\omega)}(f)=\sigma_{L^1(G)}(f), \quad
\forall f\in L^1(G,\omega)$.
\end{theorem}

\noi The last three results are based on a method developped by
Ludwig (\cite{Lu.1}). Previously, using a result of Hulanicki
(\cite{Hu.1}), Pytlik (\cite{Py.2}) had already proved the
following:

\begin{theorem}
If the weight $\omega$ satisfies
\begin{eqnarray}\label{pol}
\omega(xy) \leq C (\omega(x)+ \omega(y)), \quad \forall x,y \in G
\end{eqnarray}
for some positive constant $C$ and if $\omega^{-1}\in L^p(G)$ for
some $0<p<+\infty$, then $L^1(G,\omega)$ is symmetric.
\end{theorem}

\noi Pytlik calls a weight satisfying (\ref{pol}) a {\it
polynomial weight}. In particular, weights of the form
$$\omega(x)=(1 + \vert x \vert)^D,$$ for some positive $D$, where
$$\vert x \vert := \inf\{ n~\vert~x\in U^n\},$$ satisfy
$\omega(xy) \leq C(\omega(x)+\omega(y))$ for all $x,y \in G$, and
hence give symmetric weighted group algebras $L^1(G,\omega)$, as
Pytlik already noticed in (\cite{Py.1}).

\noi By (\cite{Py.2}) every weight satisfying $\omega(xy) \leq C
\big (\omega(x)+\omega(y)\big )$ is dominated by a weight of the
form $K(1+ \vert x \vert )^D$, $K\geq 1, D>0$. It is then easy to
check that all these weights satisfy condition (S). So the result
of Pytlik is a particular case of the result of Fendler, Gr\"
ochenig and Leinert.

\subsection{ }
Our aim is to study symmetry for weighted $L^p$-algebras. For the
rest of this section we hence assume that $(G,\omega)$ satisfies
(LPAlg), in order to be sure that $L^p(G,\omega)$ is an algebra. The
question is whether condition (S) will also imply the symmetry of
$L^p(G,\omega)$. We need some preliminary result.

\begin{lemma}\label{symmetry-o(exp)}
The weight $\omega$ satisfies condition (S) if and only if
$\omega(x)= {\mathcal O}(e^{\varepsilon \vert x\vert})$ for all
$\varepsilon>0$, where $\vert x \vert = \inf\{n~\vert~x \in
U^n\}$.
\end{lemma}

\begin{proof}
Let us assume that $\omega(x) \leq C(\varepsilon)
e^{\varepsilon\vert x\vert}$ for some constant $C(\varepsilon)$,
$\varepsilon>0$. Then

\begin{eqnarray}
\nn x\in U^k &\Rightarrow& \vert x \vert \leq k\\
\nn &\Rightarrow& \omega(x) \leq C(\varepsilon) e^{\varepsilon k}
\end{eqnarray}
and
$$\sup_{x\in U^k} \omega(x)^{\frac{1}{k}} \leq
C(\varepsilon)^{\frac{1}{k}} e^{\varepsilon}.$$ So $$1 \leq
\lim_{k\rightarrow+\infty} \sup_{x\in U^k} \omega(x)^{\frac{1}{k}}
\leq e^{\varepsilon}.$$ As this has to be true for all
$\varepsilon>0$, $\lim_{k\rightarrow+\infty} \sup_{x\in U^k}
\omega(x)^{\frac{1}{k}}=1$ and $\omega$ satisfies (S).

\noi Conversely, let us assume that there exists $\varepsilon>0$
such that $\omega(x)$ is not ${\mathcal O}(e^{\varepsilon \vert x
\vert})$. So, for every $k \in \N$, there exists $x_k \in G$ such
that $\omega(x_k)> k e^{\varepsilon\vert x_k\vert} >k$. As
$\omega$ is bounded on each $U^n$, this implies that
$"x_k\rightarrow \infty"$, which means the following: Let
$n(k):=\vert x_k\vert$. Then the sequence $(n(k))_k$ admits a
subsequence $(\tilde n(r))$ such that
\begin{eqnarray}
\nn &{}& \tilde n(r) \rightarrow +\infty, \quad \text{if }
r \rightarrow +\infty\\
\nn &{}& x_{r} \in U^{\tilde n(r)} \setminus U^{\tilde
n(r)-1} \\
\nn &{}& \omega(x_{r})> r e^{\varepsilon \tilde n(r)}> r\\
\nn &{}& \sup_{x\in U^{\tilde n(r)}} \omega(x)^{\frac{1}{\tilde
n(r)}} \geq \omega(x_{r})^{\frac{1}{\tilde n(r)}}
> r^{\frac{1}{\tilde n(r)}} e^{\varepsilon} \geq
e^{\varepsilon}.
\end{eqnarray}
Hence $$\overline {\lim}_{r\rightarrow+\infty} \sup_{x\in
U^{\tilde n(r)}} \omega(x)^{\frac{1}{\tilde n(r)}} \ge
e^{\varepsilon}>1.$$ Thus $\omega$ does not satisfy condition (S).
\end{proof}

\noi We may now prove the symmetry result for abelian groups:

\begin{theorem}\label{symmetry-lp}
Let $G$ be an abelian group such that $(G,\omega)$ satisfies
(LPAlg). Let us assume that $\omega$ satisfies condition (S). Then
$L^p(G,\omega)$ is a symmetric Banach $*$-algebra.
\end{theorem}

\begin{proof}
According to (\cite{Ku.1}), every character $\chi$ of
$L^p(G,\omega)$ is of the form $$\chi(f)=\int_G f(x)\sigma(x)dx=
\int_G \Big (f(x)\omega(x) \Big) \Big (\sigma(x) \omega^{-1}(x)
\Big )dx, \quad \forall f\in L^p(G,\omega),$$ where $\sigma$ is a
(possibly unbounded) character of the group $G$. As $f\omega \in L^p(G)$ is
arbitrary, $\vert \sigma(x) \vert \omega(x)^{-1} \in L^q(G)$, with
$\frac{1}{p} + \frac {1}{q}=1$, and, as $\omega(x) \leq
C(\varepsilon) e^{\varepsilon \vert x \vert}$ for all
$\varepsilon>0$, $$\vert \sigma(x)\vert e^{-\varepsilon \vert x
\vert} \leq \vert \sigma(x) \vert C(\varepsilon) \omega(x)^{-1}
\in L^q(G).$$ Let us assume that $\sigma$ is not unitary. Then
there exists $x_0\in G$ and $\delta>0$ such that $\vert
\sigma(x_0)\vert>1+ \delta$. By continuity, there is a non-empty open subset $V$
of $G$ such that $\vert \sigma(x)\vert > 1+ \frac{\delta}{2}$, for
all $x\in V$. Hence, for $x\in V^n$, $x=x_1 \cdot x_2 \cdots x_n$
with $x_j \in V$ for all $j$, and $\vert \sigma(x)\vert =
\prod_{j=1}^n \vert \sigma(x_j)\vert \geq (1+
\frac{\delta}{2})^n$. On the other hand, if $V\subset U^k$ where $U$ is a (relatively compact) generating neighbourhood of the identity then $V^n\subset U^{kn}$, and so $|x|\le kn$ for $x\in V^n$. This implies that, for all~$n$ (and for all $\e>0$),
\begin{eqnarray}
+\infty &>& \int_G \Big (\vert \sigma(x) \vert
e^{-\varepsilon\vert x\vert} \Big )^q dx \\
&\geq& \int_{V^n} \vert \sigma(x)\vert^q e^{-\varepsilon q\vert
x\vert}dx\\
&\geq& \big(1+\frac{\delta}{2}\big)^{qn} \int_{V^n} e^{-\varepsilon q
kn}dx\\
&\geq& \Big( \big(1+\frac{\delta}{2}\big) \cdot e^{-\varepsilon k} \Big)^{qn} |V|. \label{char_estimate}
\end{eqnarray}
But we can choose $\e$ so that $\big(1+\frac{\delta}{2}\big) \cdot
e^{-\varepsilon k} >1$, then \eqref{char_estimate} tends to
$+\infty$ with $n$. This is a contradiction which shows that
$\sigma$, and hence $\chi$ are unitary. So $L^p(G,\omega)$ is
symmetric.
\end{proof}

\noi {\bf Example:} $L^2(\R, e^{\sqrt{\vert x\vert}})$ is a
symmetric Banach $*$-algebra (section \ref{subexp-w}).

\subsection{ }
Before studying the non-abelian case, let us first recall the
generalized Minkowski relation: Let $X,Y$ be measure spaces and
let $F$ be a measurable function on $X \times Y$. Then, for all
$p\geq 1$, $$\Bigg ( \int_X \Big ( \int_Y \vert F(x,y) \vert dy
\Big )^pdx \Bigg )^{\frac{1}{p}} \leq \int_Y \Big ( \int_X \vert
F(x,y) \vert^p dx \Big )^{\frac{1}{p}} dy.$$ We need the following
relation:

\begin{lemma}\label{norminequality}
Let us assume that the weight $\omega$ satisfies (LPAlg) and is polynomial in the sense
of Pytlik, i. e. that it satisfies $$\omega(xy) \leq C \Big (
\omega(x) + \omega(y) \Big ), \quad \forall x,y \in G,$$ for some
constant $C>0$. Then $$\Vert f*g \Vert_{p,\omega} \leq C \Big (
\Vert f \Vert_{p,\omega}~\Vert g\Vert_1 + \Vert g
\Vert_{p,\omega}~\Vert f\Vert_1 \Big), \quad \forall f,g \in
L^p(G,\omega) \subset L^1(G).$$
\end{lemma}

\begin{proof}
\begin{eqnarray}
\nn \Vert f*g\Vert_{p,\omega} &=& \Big ( \int_G \vert \int_G f(y)
g(y^{-1}x)dy \vert^p \omega(x)^p dx \Big )^{\frac{1}{p}}\\
\nn &\leq& C \Big (\int_G \vert \int_G f(y) g(y^{-1}x) \big
(\omega(y)
+ \omega(y^{-1}x) \big ) dy \vert^p dx \Big )^{\frac{1}{p}}\\
\nn &\leq& C \Big (\int_G \vert \int_G f(y) g(y^{-1}x) \omega(y)dy
\vert^p dx \Big )^{\frac{1}{p}} + C \Big ( \int_G \vert \int_G
f(y) g(y^{-1}x) \omega(y^{-1}x) dy \vert^p dx \Big
)^{\frac{1}{p}}\\
\nn &=& I + II
\end{eqnarray}
by the triangle inequality for $\Vert \cdot \Vert_{p}$. Then the
generalized Minkowski inequality implies that
\begin{eqnarray}
\nn I &\leq& C \Big ( \int_G \big ( \int_G \vert f(y)
g(y^{-1}x)\vert \omega(y)dy \big )^p dx \Big)^{\frac{1}{p}} \\
\nn &\leq& C \Big ( \int_G \big (\int_G \vert f(xu^{-1}) g(u)
\vert \omega(xu^{-1}) du \big )^pdx \Big )^{\frac{1}{p}} \quad
(y=xu^{-1})\\
\nn &\leq& C \int_G \Big (\int_G \vert f(xu^{-1}) \vert^p \vert
g(u)\vert^p \omega(xu^{-1})^pdx \Big)^{\frac{1}{p}}du \\
\nn &=& C \Vert f\Vert_{p,\omega}~\Vert g\Vert_1
\end{eqnarray}
and
\begin{eqnarray}
\nn II &\leq& C \Big ( \int_G \big ( \int_G \vert f(y) g(y^{-1}x)
\vert \omega(y^{-1}x) dy \big )^p dx \Big )^{\frac{1}{p}}\\
\nn &\leq& C \int_G \Big ( \int_G \vert f(y) \vert^p \vert
g(y^{-1}x) \vert^p \omega(y^{-1}x)^pdx \Big )^{\frac{1}{p}} dy \\
\nn &=& C \Vert g \Vert_{p,\omega} ~\Vert f \Vert_1.
\end{eqnarray}
\end{proof}
\noi We may now use the methods of Pytlik (\cite{Py.1},
\cite{Py.2}) to show the symmetry of the algebra $L^p(G,\omega)$
for polynomial weights. This is done via the following lemmas.

\begin{lemma}
Let $(G,\omega)$ satisfy (LPAlg). Then, for any $f\in L^p(G,\omega)
\subset L^1(G)$, $r_1(f) \leq r_{p,\omega}(f)$, where $r_1(f)$
denotes the spectral radius of $f$ in $L^1(G)$ and
$r_{p,\omega}(f)$ denotes the spectral radius of $f$ in
$L^p(G,\omega)$.
\end{lemma}
\begin{proof}
From $$\Vert f\Vert_1 = \int_G \vert f(x) \vert \omega(x)
\omega^{-1}(x) dx \leq \big ( \int_G \omega^{-q}(x)dx \big
)^{\frac{1}{q}} \Vert f \Vert_{p,\omega} = C \Vert f
\Vert_{p,\omega}$$ we deduce $$\Vert f^{*n}\Vert_1^{\frac{1}{n}}
\leq C^{\frac{1}{n}} \Vert f^{*n}
\Vert_{p,\omega}^{\frac{1}{n}}~,$$ where $f^{*n} = f*f* \cdots *f$
($n$ factors). Hence, for $n\rightarrow +\infty$, $$r_1(f) \leq
r_{p,\omega}(f).$$
\end{proof}

\begin{lemma}
Let $(G,\omega)$ satisfy (LPAlg). Let us assume that $\omega$ is
polynomial in the sense of Pytlik, i. e. that $$\omega(xy) \leq C
\big (\omega(x) + \omega(y) \big), \quad \forall x,y \in G.$$ Then
$$r_1(f)=r_{p,\omega}(f), \quad \forall f\in L^p(G,\omega).$$
\end{lemma}

\begin{proof}
By the methods of Pytlik (\cite{Py.2}), (\ref{norminequality}) gives
$$\Vert f*f\Vert_{p,\omega} \leq 2C
\Vert f \Vert_{p,\omega}~\Vert f\Vert_1$$ and, by induction,
$$\Vert f^{*2^n} \Vert_{p,\omega} \leq (2C)^n \Vert f
\Vert_{p,\omega}~\Vert f\Vert_1^{2^n-1}.$$ So,
\begin{eqnarray}
\nn r_{p,\omega}(f) &=& \lim_{n\rightarrow +\infty} \Vert
f^{*2^n}\Vert_{p,\omega}^{2^{-n}}\\
\nn &\leq& \lim_{n\rightarrow +\infty} (2C)^{n \cdot 2^{-n}} \Vert
f
\Vert_{p,\omega}^{2^{-n}}~\Vert f\Vert_1^{1-2^{-n}}\\
\nn &=& \Vert f \Vert_1.
\end{eqnarray}
Finally, $$r_{p,\omega}(f)= r_{p,\omega}(f^{*n})^{\frac{1}{n}}
\leq \Vert f^{*n}\Vert_1^{\frac{1}{n}}, \quad \forall n,$$ and
$$r_{p,\omega}(f) \leq \lim_{n\rightarrow +\infty} \Vert f^{*n}
\Vert_1^{\frac{1}{n}} = r_1(f).$$
\end{proof}

\noi We finally get:

\begin{theorem}
Let $(G,\omega)$ satisfy (LPAlg). Let us assume that $\omega$ is
polynomial in the sense of Pytlik, i. e. that $$\omega(xy) \leq C
\big ( \omega(x) + \omega(y) \big ), \quad \forall x,y \in G.$$
Then $L^p(G,\omega)$ is a symmetric Banach $*$-algebra.
\end{theorem}

\begin{proof}
This follows from Lemma 3.1 of (\cite{Fe.Groe.Lei.}) applied to
${\mathcal A}:= L^p(G,\omega)$ and ${\mathcal B}:= L^1(G)$, and
from the result of Losert (\cite{Lo.}) about the symmetry of
$L^1(G)$. As a matter of fact, these results imply that
$$\sigma_{L^1(G)}(f)=\sigma_{L^p(G,\omega)}(f), \quad \forall
f=f^* \in L^p(G,\omega),$$ and, as $L^1(G)$ is symmetric,
$\sigma_{L^1(G)}(f) \subset \R$.
\end{proof}

\noi {\bf Example:} We know that for $D>0$ sufficiently large,
$\omega(x)=(1+ \vert x\vert)^D$, where $\vert
x\vert=\inf\{n~\vert~x\in U^n\}$, gives rise to an $L^p$-algebra
(\ref{poly-w}). Hence this algebra is symmetric.

\medskip
\noi One may conjecture that, more generally, $L^p(G,\omega)$ is a
symmetric Banach $*$-algebra, if $(G,\omega)$ satisfies (LPAlg) and
$\omega$ satisfies condition (S), or, equivalently, the
GRS-condition. But this remains an open question.

\section{Functional calculus}
\subsection{}
Let $(G,\omega)$ satisfy (LPAlg). The aim of the following section
is the construction of functional calculus for all continuous
functions $f$ with compact support such that $f=f^*$, where
$f^*(x)= \overline {f(x^{-1})}$, for all $x\in G$. We will follow
the method of (\cite{Hu.2}, \cite{Dz.Lu.Mo.}) and use their results. To use this method, we have to
bound $$u(nf):= \sum_{k=1}^{\infty} \frac{i^k}{k!} n^k f^{*k}$$ in
$L^p(G,\omega)$ and show that there are "enough" functions
$\varphi: \R \rightarrow\R$, periodic with period $2\pi$, with
$\varphi(0)=0$, such that $$\varphi \{f\} := \sum_{n\in \Z} u(nf)
\hat \varphi(n)$$ converges in $L^p(G,\omega)$. For more details
on functional calculus, see among others (\cite{Di.},
\cite{Dz.Lu.Mo.}).

\noi Let us first recall that, for all continuous functions $g$
with compact support, $\Vert g \Vert_1 \leq \Vert g
\Vert_{1,\omega}$ and $\Vert g \Vert_1 \leq C \Vert g
\Vert_{p,\omega}$ for some positive constant $C$. Moreover,
$$\Vert u(nf)\Vert \leq \sum_{k=1}^{+\infty} \frac{1}{k!} n^k
\Vert f\Vert^k \leq e^{n\Vert f\Vert}$$ in any convenient norm
$\Vert \cdot \Vert$. So for any continuous function $f$ with
compact support such that $f=f^*$, the series defining $u(nf)$
converges in $L^1(G)$, $L^1(G,\omega)$ and $L^p(G,\omega)$ to the
same element, i. e. the notation $u(nf)$ represents a function
belonging to $L^1(G,\omega) \cap L^p(G,\omega) \subset L^1(G)$. We
will now deduce a bound for $\Vert u(nf) \Vert_{p,\omega}$ from
the bound for $\Vert u(nf)\Vert_{1,\omega}$ which was established
in (\cite{Dz.Lu.Mo.}). Let us recall the following notations and
facts from (\cite{Dz.Lu.Mo.}): There exists a constant $C>1$ such
that $$\omega(x) \leq e^{C \vert x\vert},~\forall x\in G,$$ where
$\vert x\vert= \inf \{n~\vert~x\in U^n\}$ for an arbitrary
(relatively compact) generating neighbourhood $U$. We denote
\begin{eqnarray}
\nn s(n)&:=& \sup_{x\in U^n} \omega(x),~\quad \forall
n\in \N^*\\
\nn s(0) &:=& 1\\
\nn \omega_1(x) &:=& s(\vert x\vert)\\
\nn \omega_2(x) &:=& e^{C \vert x\vert} \\
\nn s_2(n) &:=& \sup_{x\in U^n} \omega_2(x)=e^{Cn}.
\end{eqnarray}
We also consider an arbitrary increasing function $r:\N\rightarrow
\N$, which will be specified later. It is shown in
(\cite{Dz.Lu.Mo.}) that there exist positive constants $C_1, C_2$
such that
\begin{equation}\label{unf-bound}
\Vert u(nf)\Vert_{1,\omega} \leq C_1 (1+ \vert n \vert)(1+
\vert n \vert r(\vert n \vert))^{\frac{Q}{2}} s(\vert n \vert
r(\vert n \vert)) e^{C_2 \big ( \frac{\vert n \vert}{s_2(r(\vert n
\vert))}\big)}, \quad \forall n \in \Z,
\end{equation}
where $Q$ denotes the power appearing in the polynomial growth
condition of the group $G$, i. e. $\vert U^n \vert \leq K n^Q$,
for all $n\in \N^*$, for some positive constant $K$.

\noi From the formal representation $u(f)=e^{if}-1$ we get the
following identity valid also in the non-unital case:
\begin{align*}
u(nf) &= e^{inf}-1 = e^{i(n-1)f}*e^{if}-1 = (u((n-1)f)+1)*(u(f)+1)-1=\\
&=u((n-1)f)*u(f)+u((n-1)f)+u(f).
\end{align*}
By induction it follows that
$$
u(nf) = nu(f) + \sum_{k=1}^{n-1} u(kf)*u(f).
$$
From \eqref{Lp-L1-module}, we get an estimate
$$
\|u(nf)\|_{p,\omega} \le n\|u(f)\|_{p,\omega} + \sum_{k=1}^{n-1} \|u(kf)\|_{1,\omega} \|u(f)\|_{p,\omega}.
$$

\noi Using the bound for $\Vert u(kf)\Vert_{1,\omega}$ obtained in
(\cite{Dz.Lu.Mo.}) and recalled in \eqref{unf-bound}, we get
\begin{eqnarray}
\nn \Vert u(nf)\Vert_{p,\omega} &\leq& Kn + K_1
\sum_{k=1}^{n-1}(1+k)(1+kr(k))^{\frac{Q}{2}} s\big ( kr(k)\big )
e^{C_2 \big ( \frac{k}{s_2(r(k))}\big )} \\
\nn &\leq& Kn + K_1 (n+1)(1+nr(n))^{\frac{Q}{2}} s\big (nr(n)\big
) \sum_{k=1}^{n-1} e^{C_2 \big ( \frac{k}{s_2(r(k))} \big )},
\end{eqnarray}
for some constants $K$, $K_1$, as the functions $s$ and $r$
are increasing. (Here $1+k$ could be bounded by $n$ as well, but we choose $n+1$ to comply with assumptions of \cite{Dz.Lu.Mo.}). As in ({\cite{Dz.Lu.Mo.}), we put $r(n):= \ln (\ln
n)+1$, for $n\geq e^e$. Hence, for $k\geq \max(e^e, e^C)$,
$$s_2(r(k))= e^{Cr(k)} \geq e^{C \ln (\ln k)} = (\ln k)^C$$ and
$$e^{C_2 \big ( \frac{k}{s_2(r(k))}\big )} \leq e^{C_2
\frac{k}{(\ln k)^C}} \leq e^{C_2 \frac{n}{(\ln n)^C}},$$ as the
function $f(x)= \frac{x}{(\ln x)^C}$ is increasing for $x\geq
e^C$. This allows to estimate the sum over $k$ by $ne^{\frac{C_2n}{(\ln n)^C}}$.
Moreover, as in (\cite{Dz.Lu.Mo.}), for $n\geq e^e$,
\begin{eqnarray}
\nn &{}&\ln (\ln n) \leq r(n) \leq 2 \ln (\ln n) \leq 2n\\
\nn &{}&s\big (nr(n)\big ) \leq s(n)^{r(n)} \leq s(n)^{2\ln (\ln
n)}.
\end{eqnarray}
Finally, noticing that $nf=(-n)(-f)$, we may compute $\Vert
u(nf)\Vert_{p,\omega}$ even for negative $n$ (by replacing the
constants depending on $f$ by the sup of the corresponding
constants for $f$ and $-f$). So, there exist positive constants
$A_1, A_2$ (depending on $f$ and $\omega$) such that $$\Vert
u(nf)\Vert_{p,\omega} \leq A_1 (1+|n|)^2 (1+n^2)^{\frac{Q}{2}} s(\vert n\vert)^{2\ln (\ln \vert n \vert)}
e^{A_2 \big ( \frac{\vert n\vert}{(\ln \vert n\vert )^C}\big )}$$
for all $\vert n\vert \geq \max(e^e, e^C)$. We thus obtain a
similar bound as for $\Vert u(nf)\Vert_{1,\omega}$, except that
the factor $(1+\vert n\vert)$ has been replaced by $(1+|n|)^2$. Of course the constants are slightly different too.
They depend on $f$ and $\omega$. We may conclude exactly as in
(\cite{Dz.Lu.Mo.}).

\subsection{}
Let us recall the non-abelian Beurling-Domar condition (BDna)
given by $$\sum_{n\in \N, n\geq e^e} \frac{\big (\ln (\ln n) \big
) \ln \big ( s(n)\big )}{1+n^2} <+\infty.$$ It is independent of
the choice of the generating neighbourhood $U$ used to compute
$s(n)$. See (\cite{Dz.Lu.Mo.}) for more details on that condition.
We then have the following result:

\begin{theorem}\label{calculus}
Let $(G,\omega)$ satisfy (LPAlg). Let us assume that moreover the
weight $\omega$ satisfies the (BDna) condition. Let $f=f^*$ be a
continuous function with compact support. Then, given
$a,b,\varepsilon$ such that $0<a<a+\varepsilon<
b-\varepsilon<b<2\pi$, there exists a function $\psi:\R
\rightarrow \R$, continuous, periodic of period $2\pi$ such that
$\text{supp}\psi \cap [0,2\pi] \subset [a,b]$, $\psi \equiv 1$ on
$[a+\varepsilon, b-\varepsilon]$ and $$\sum_{n\in \Z} \Vert
u(nf)\Vert_{p,\omega} \vert \hat \psi(n)\vert < +\infty.$$ Hence
this defines a function $$\psi\{f\}:= \sum_{n\in \Z} \hat \psi(n)
u(nf) \in L^p(G,\omega) \cap L^1(G,\omega)$$ and the properties of
functional calculus are satisfied, i. e. $$\chi(\psi\{f\}) =
\psi(\chi(f))$$ for every character $\chi$ of the abelian Banach
$*$-subalgebra of $L^p(G,\omega)$ generated by $f$,
\begin{eqnarray} \nn &{}&\pi
(\psi\{f\}) = \psi(\pi(f)), \quad \forall \pi \in
\widehat{L^p(G,\omega)} \equiv \hat G,\\
\nn &{}&(\varphi\psi)\{f\} = \varphi\{f\} * \psi\{f\},
\end{eqnarray}
if the functions $\varphi$ and $\varphi\psi$ still have the
correct properties to allow functional calculus.
\end{theorem}
\begin{proof}
See (\cite{Dz.Lu.Mo.}), pages 337 to 345. Here we use again the
argument that if a series converges in $L^1(G,\omega)$ and
$L^p(G,\omega)$, then it also converges in $L^1(G)$ and the limit
is the same in the three spaces.
\end{proof}

\subsection{}
{\bf Examples:} a) If $G \in [PG]$ and $\omega(x)=K(1+\vert
x\vert)^D$ for $K \geq 1$ and $D>0$ large enough, then $(G,
\omega)$ satisfies (LPAlg) (\ref{poly-w}). It is easy to check, that
$\omega$ also verifies (BDna) and so functional calculus exists.

\noi b) Let $(G,\omega)$ satify (LPAlg) and let us assume that
$\omega(xy)\leq C(\omega(x)+\omega(y))$ for all $x,y\in G$. By
(\cite{Py.2}), such a weight is bounded by a weight of the form
$K(1+\vert x\vert)^D$ and hence (BDna) is verified. Functional
calculus exists.

\noi c) If $G\in [PG]$, then $$\omega(x):= e^{C\vert
x\vert^{\gamma}}, \quad 0<\gamma<1,$$ is such that $(G,\omega)$
satisfies (LPAlg) (\ref{subexp-w}) and the weight $\omega$ verifies (BDna)
(\cite{Dz.Lu.Mo.}). Functional calculus exists.

\medskip
\noi {\bf Remarks:} a) Condition (BDna) is independent of the
choice of the generating neighbourhood $U$.

\noi b) Condition (BDna) is only slightly more restrictive than
the well known Beurling-Domar condition in the abelian case.

\noi c) If $\omega$ satisfies (BDna), it also verifies condition
(S).

\medskip \noi See (\cite{Dz.Lu.Mo.}) for more details.

\medskip \noi Functional calculus is a very useful tool to prove
different harmonic analysis properties, as will be shown in the
rest of this paper.

\section{Regularity}

\subsection{}
For abelian Banach algebras, regularity is defined as follows by \v{S}ilov \cite{shilov}: Let
${\mathcal A}$ be an abelian Banach algebra and let
$\Delta({\mathcal A})$ denote the space of characters of ${\mathcal
A}$. Then ${\mathcal A}$ is said to be regular if, given any
$\varphi \in \Delta({\mathcal A})$ and any closed set $F \subset
\Delta({\mathcal A})$ not containing $\varphi$, there exists $x\in
{\mathcal A}$ such that $\hat x(\varphi) = \varphi(x)=1$ and $\hat
x\vert_F \equiv 0$, where $\hat x$ denotes the Gelfand transform
of $x$.

\noi In the non-abelian case $\Delta(\cal A)$ should be replaced
by the space $Prim_*\cal A$ defined as the set of all kernels of
topologically irreducible $*$-representations of $\cal A$. The set
$Prim_*\cal A$ is equipped with the hull-kernel topology.

\subsection{}
As previously, we assume that $(G,\omega)$ satisfies (LPAlg). Let
$C_c(G)$ denote the set of continuous functions with compact
support on $G$. It is obvious that $C_c(G) \subset L^p(G,\omega)
\subset L^1(G)$, that $C_c(G)$ is dense in $L^p(G,\omega)$ and in $L^1(G)$, that for all $\pi \in \hat G\equiv
\widehat{L^p(G,\omega)}$,
$$\Vert \pi(f)\Vert_{op} \leq \Vert f\Vert_1 \leq C \Vert f
\Vert_{p,\omega}, \quad \forall f\in L^p(G,\omega),$$ resp. $\Vert
\pi(f)\Vert_{op} \leq \Vert f \Vert_1$ for all $f\in L^1(G)$.
Moreover, we have seen in the previous section that functional
calculus is possible on the self-adjoint elements of $C_c(G)$,
provided the weight $\omega$ satisfies (BDna). This imlies that
the arguments of (\cite{Dz.Lu.Mo.}, pages 350 and 351) remain
valid. In particular, we have the following results:

\begin{theorem}
Let $(G,\omega)$ satisfy (LPAlg). Let us assume that the weight
$\omega$ verifies (BDna). We then have:

\noi (i) The map \begin{eqnarray} \nn \Psi:
\text{Prim}_*L^1(G)&\rightarrow& \text{Prim}_* L^p(G,\omega)\\
\nn \text{ker} \pi &\mapsto& \text{ker}\pi \cap L^p(G,\omega)
\end{eqnarray}
is a homeomorphism.

\noi (ii) In particular, $\text{Prim}_*L^p(G,\omega)$,
$\text{Prim}_*L^1(G,\omega)$ and $\text{Prim}_*L^1(G)$ are
homeomorphic.

\noi (iii) Given any $\rho \in \hat G$ and any open neighbourhood
$N$ of $\rho$ in $\hat G$, resp. any open neighbourhood $N_1$ of
$\text{ker}\rho \cap L^p(G,\omega)$ in
$\text{Prim}_*L^p(G,\omega)$, there exists $f \in L^p(G,\omega)$
such that $\rho(f)\neq 0$ and $\pi(f)=0$ for all $\pi \in \hat
G\setminus N$, resp. for all $\pi$ such that $\text{ker}\pi \cap
L^p(G,\omega) \in \text{Prim}_* L^p(G,\omega)\setminus N_1$.
\end{theorem}
\begin{proof} See (\cite{Dz.Lu.Mo.}). Part of the argument relies
heavily on functional calculus and on the $*$-regularity of groups
with polynomial growth.
\end{proof}

\noi Point (iii) of the previous theorem, which is often called
Domar's property, corresponds to the regularity of abelian Banach
algebras.

\section{Weak Wiener property}
Let us recall the following definitions:

\begin{definition}
Let ${\mathcal A}$ be a Banach algebra.

\noi (i) A representation $(T,V)$ of ${\mathcal A}$ on a vector
space $V$ is said to be algebraically irreducible, if there are no
non-trivial $T$-invariant subspaces in $V$.

\noi (ii) The algebra ${\mathcal A}$ is said to have the weak
Wiener property, if every proper closed two-sided ideal of
${\mathcal A}$ is contained in the kernel of an algebraically
irreducible representation.
\end{definition}

\noi Let $(G,\omega)$ satisfy (LPAlg) and let $(f_s)_s$ be an
approximate unit of $L^p(G,\omega)$ with the properties discussed
in \ref{sec-approx}.

\noi In (\cite{Dz.Lu.Mo.}) it is shown that, provided $\omega$
satisfies (BDna), there exists a periodic function $\varphi$ of
period $2\pi$ with $\varphi(1)=1$, $\varphi \equiv 0$ in a
neighbourhood of $0$, such that $\varphi\{f_s\}$ is defined in
$L^1(G,\omega)$. By our section on functional calculus in
$L^p(G,\omega)$, the same $\varphi\{f_s\}$ also converges in
$L^p(G,\omega)$, i. e. $\varphi\{f_s\} \in L^1(G,\omega) \cap
L^p(G,\omega)$ for all $s$. Moreover, in (\cite{Dz.Lu.Mo.}) it is
shown that $$\Vert \varphi\{f_s\} *f - f \Vert_{1,\omega}
\rightarrow 0$$ for all continuous functions $f$ with compact
support in $G$. Hence, for any $f,g \in C_c(G)$, $f*g \in C_c(G)
\subset L^1(G,\omega) \cap L^p(G,\omega)$ and $$\Vert
\varphi\{f_s\}*f*g -f*g\Vert_{p,\omega} \leq \Vert
\varphi\{f_s\}*f-f\Vert_{1,\omega} \Vert g\Vert_{p,\omega}.$$ So
\begin{eqnarray}\label{limit}
\Vert \varphi\{f_s\}*f*g -f*g \Vert_{p,\omega} \rightarrow 0.
\end{eqnarray}
This gives the following result:

\begin{lemma}
Under the assumptions above, let $I$ be a proper closed two-sided ideal of $L^p(G,\omega)$.
Then there exists $s$ such that $\varphi\{f_s\} \notin I$.
\end{lemma}

\begin{proof} Let us assume that $\varphi\{f_s\} \in I$, for all
$s$. Then $\varphi\{f_s\}*f*g \in I$ for all $s$ and all $f,g\in
C_c(G)$. As $I$ is closed, the relation (\ref{limit}) shows that
$f*g\in I$ for all $f,g\in C_c(G)$. But this implies that
$I=L^p(G,\omega)$, by density, which is a contradiction.
\end{proof}
\noi We are now able to prove the weak Wiener property:

\begin{theorem}
Let $(G,\omega)$ be (LPAlg). Let us also assume that the weight
$\omega$ satisfies (BDna). Then the algebra $L^p(G,\omega)$ has
the weak Wiener property.
\end{theorem}
\begin{proof} The proof is standard, but we repeat it for the sake
of completeness. Let $I$ be a proper, closed, two-sided ideal of
$L^p(G, \omega)$. Let $s$ and $\varphi$ be such that
$\varphi\{f_s\} \notin I$. Let $\psi$ be another function such
that functional calculus $\psi\{f_s\}$ is possible and such that
$\psi\equiv 1$ on the support of $\varphi$. Such a $\psi$ exists
by theorem \ref{calculus}. Then $$\psi\{f_s\}*\varphi\{f_s\} =
(\psi\varphi)\{f_s\} = \varphi\{f_s\}.$$ Let us consider the
algebra ${\mathcal A}:= L^p(G,\omega)/I$. We have
\begin{eqnarray} \nn &{}& 0\neq \dot{\wideparen{\varphi\{f_s\}}} \in {\mathcal
A}\\
\nn &{}& (\dot{\wideparen{\psi\{f_s\}}}-1)* \dot{\wideparen{\varphi\{f_s\}}}=0 \quad
\text{in } {\mathcal A} \oplus \C,
\end{eqnarray}
where the dot denotes the equivalence class in the quotient space
$L^p(G,\omega)/I$. Hence, by (\cite{Bo.Du.})
$\dot{\wideparen{\psi\{f_s\}}} -1$ is not invertible in ${\mathcal
A} \oplus \C$, i. e. $1\in \sigma_{\mathcal A}
(\dot{\wideparen{\psi\{f_s\}}})$. So
$\dot{\wideparen{\psi\{f_s\}}} \notin \text{rad}({\mathcal A})$,
where $\text{rad}({\mathcal A})$ denotes the radical of ${\mathcal
A}$. This implies that there exists an algebraically irreducible
representation $(\tilde T,V)$ of ${\mathcal A}$ such that $\tilde
T(\dot{\wideparen{\psi\{f_s\}}})\neq 0$. We then define the
non-trivial algebraically irreducible representation $(T,V)$ of
$L^p(G,\omega)$ by $T(f):=\tilde T(\dot {f})$. By construction, $I
\subset \text{ker}T$. Hence $L^p(G,\omega)$ is weakly Wiener.
\end{proof}

\section{Wiener property}
We start with the following definition:

\begin{definition}
Let ${\mathcal A}$ be a Banach $*$-algebra.

\noi (i) A representation $(T,V)$ of ${\mathcal A}$ on a Banach
space $V$ is said to be topologically irreducible, if there are no
non-trivial closed $T$-invariant subspaces in $V$. In particular, this
definition is applied to unitary representations $(T, {\mathcal
H})$ on Hilbert spaces ${\mathcal H}$.

\noi (ii) The algebra ${\mathcal A}$ is said to have the Wiener
property, if every proper closed two-sided ideal in ${\mathcal A}$ is
contained in the kernel of a topologically irreducible unitary
representation of ${\mathcal A}$.
\end{definition}

\noi It is well known that every symmetric Banach $*$-algebra
which has the weak Wiener property, also has the Wiener property
(see \cite{Le.1} and \cite{Le.2}). This leads us to the following
result:

\begin{theorem} Let $(G,\omega)$ satisfy (LPAlg). Then the algebra
$L^p(G,\omega)$ has the Wiener property in the following cases:

\noi a) $G$ is abelian and $\omega$ satisfies $(BDna)$.

\noi b) $G$ is non abelian and $\omega(xy) \leq
C(\omega(x)+\omega(y))$, for all $x,y\in G$, for some constant
$C>0$.

\noi c) $G$ is non abelian and $\omega(x)=K(1+\vert x\vert)^D$ for
$D$ large enough, with $\vert x\vert=\inf\{n~\vert~x\in U^n\}$,
where $U$ is an arbitrary, relatively compact, generating
neighbourhood of the identity.
\end{theorem}

\medskip \noi If one could prove that property (S) implies the
symmetry of the algebra $L^p(G,\omega)$ (conjecture), then the
property (BDna) would imply the Wiener property.

\section{Minimal ideals of a given hull}

As always, we assume that $(G,\omega)$ satisfies (LPAlg). We also
suppose that $\omega$ satisfies (BDna). For details of the
following we refer to (\cite{Dz.Lu.Mo.}). We will use the
following notations:

\noi Let us denote by $\Phi$ the set of functions $\varphi$ from
$\R$ to $\R$, periodic of period $2\pi$, with $\varphi(0)=0$, with
$\text{supp}\varphi \cap [0,2\pi]$ compact contained in
$]0,2\pi[$, which operate on the set of continous, self-adjoint
functions with compact support in the algebras $L^1(G,\omega)$ and
$L^p(G,\omega)$. The construction of such functions is described
in more details in (\cite{Dz.Lu.Mo.}). See also theorem
\ref{calculus}.

\noi For any closed ideal $I$ of $L^p(G,\omega)$, we define the
hull of $I$ by $$h(I):= \{ \text{ker}\pi \in \text{Prim}_*
L^p(G,\omega)~\vert~I\subset \text{ker}\pi\}.$$ For any compact
subset $C$ of $\text{Prim}_*L^p(G,\omega)$ (endowed with the
hull-kernel topology), let us define
\begin{eqnarray}
\nn \tilde C &:=& \{\pi \in \hat G~\vert~\text{ker}\pi \in C\}
\quad \text{ and } \quad \Vert f\Vert_C := \sup_{\pi \in \tilde C}
\Vert \pi (f)
\Vert_{op}\\
\nn m(C) &:=& \{ \varphi\{f\}~\vert~f=f^*, f\in C_c(G), \Vert f
\Vert_1\leq 1, \varphi\in \Phi, \varphi\equiv 0 \text{ on a
neighbourhood of } [-\Vert f\Vert_C, \Vert f\Vert_C]\}
\end{eqnarray}
Let $j(C)$ be the closed two-sided ideal of $L^p(G,\omega)$
generated by $m(C)$. As in (\cite{Dz.Lu.Mo.}), one may prove:

\begin{lemma} The hull of $j(C)$ is $C$.
\end{lemma}

\begin{proof} See (\cite{Dz.Lu.Mo.}).
\end{proof}

\begin{theorem}
Let $(G,\omega)$ satisfy (LPAlg). Let the weight $\omega$ satisfy
(BDna). We also assume that the algebra $L^p(G,\omega)$ is
symmetric. Let $C$ be a closed subset of
$\text{Prim}_*L^p(G,\omega)$ ($\equiv \text{Prim}_*L^1(G,\omega)
\equiv \text{Prim}_*L^1(G)$). There exists a closed two-sided
ideal $j(C)$ of $L^p(G,\omega)$ with $h(j(C))=C$, which is
contained in every two-sided closed ideal $I$ with $h(I)=C$.

\noi This is in particular the case if $G$ is either abelian and
$\omega$ satisfies (BDna) or if $G$ is non-abelian and $\omega$
satisfies the relationship $\omega(xy) \leq C \big
(\omega(x)+\omega(y)\big)$ for all $x,y \in G$, for some positive
constant $C$.
\end{theorem}

\begin{proof}
See (\cite{Dz.Lu.Mo.}).
\end{proof}

\section{A symmetric algebra having infinite-dimensional irreducible representations}

Using C. Read's example of a quasi-nilpotent operator on $\ell_1$
with no nontrivial invariant subspaces, one can construct an
example of a weighted algebra $L^p(\R,\omega)$ for any $p>1$ which
is symmetric and has infinite-dimensional topologically
irreducible representations. Up to our knowledge, the first
application of this type is due A.~Atzmon (\cite{atzmon}).

\noi Let $T:\ell_1\to\ell_1$ be the operator constructed in \cite{read-nilpotent}. It is quasi-nilpotent, so that
\begin{equation}\label{estimate-T}
\|T\|\le1,\qquad
\|T^{n}\|^{1/n}\to 0, \quad n\to\infty.
\end{equation}

\subsection{Definition of the weights and symmetry}

We will introduce a family of weights on $\R$ which will depend on $p\ge1$. For $p=1$, we take as a weight
\begin{equation}\label{w-T-nilpotent}
\omega(x) = \max\{\, \|e^{xT}\|, \,\|e^{-xT}\| \}.
\end{equation}

\noi Obviously, $\omega$ is submultiplicative. Thus,
$L^1(\R,\omega)$ is an algebra.
For $p>1$, we put $\omega_1(x) = \omega(x) (1+|x|)^2$; 
by sections \ref{poly-w} and \ref{product-w}, this is an $L^p$-algebra weight (possibly after multiplication by a constant).
\medskip

\noi Next we prove the following estimate:

\blm Let $\omega$ be defined by \eqref{w-T-nilpotent}. Then
$\omega(x)=\mathcal O(\exp(\e |x|))$, $x\to\infty$, for any
$\e>0$. \elm \bpr Let $\e>0$ be given. It is enough to estimate
$\omega(x)$ for $x>0$. By \eqref{estimate-T}, there is $N(\e)$
such that $\|T^n\|<\e^n$ for all $n\geq N(\e)$. Separate the
series into two parts: $\omega(x) = \Omega_1(x) + \Omega_2(x)$,
where
\begin{align*}
\Omega_1(x) &= \sum_{n< N(\e)} {\|T^n\|x^n\over n!}, \\
\Omega_2(x) &= \sum_{n\ge N(\e)} {\|T^n\|x^n\over n!}.
\end{align*}
For $\Omega_2$, we have:
$$
\Omega_2(x) \le \sum_{n\ge N(\e)} {\e^n x^n\over n!} \le \exp(\e x).
$$
For $\Omega_1$:
$$
\Omega_1(x) \le  \sum_{n< N(\e)} {x^n\e^n\over n!\e^n}\le \e^{-N(\e)}\sum_{n< N(\e)} {x^n\e^n\over n!} \le \e^{-N(\e)}\exp(\e x).
$$
Now,
$$
\omega(x) \le \exp(\e x)\, \big( 1+\e^{-N(\e)}\big) = C(\e)\exp(\e x),
$$
what proves the lemma.
\epr

\bcor The algebra $L^1(\R,\omega)$ and every algebra
$L^p(\R,\omega_1)$ with $p>1$ are symmetric. \ecor \bpr It is easy
to see that $\omega_1(x)=\mathcal O(\exp(\e x))$ for any $\e>0$ as
well. By lemma \ref{symmetry-o(exp)}, $\omega$ and $\omega_1$
satisfy condition (S), so by theorems \ref{symmetry-l1} and
\ref{symmetry-lp} the algebras $L^1(\R,\omega)$ and
$L^p(\R,\omega_1)$ are symmetric. \epr

\subsection{An infinite-dimensional irreducible representation}

Let ${\mathcal A}$ stand for $L^1(\R,\omega)$ if $p=1$ and for
$L^p(\R,\omega_1)$ if $p>1$. Now we can put
$$
U(f) = \int_\R \exp(xT)f(x)dx
$$
for any $f\in {\mathcal A}$. First of all, we will show that this
integral converges absolutely. We can estimate
$$
\|U(f)\| \le \int_\R \|\exp(Tx)\|\, |f(x)| dx \leq \int_\R
\omega(x) |f(x)| dx.
$$
If $p=1$, this equals $\|f\|_{1,\omega}$. If $p>1$,
$$
\int_\R \omega(x) |f(x)| dx = \int_\R \omega_1(x) |f(x)| (1+|x|)^{-2} dx \le \|\omega_1 f\|_p \|(1+|x|)^{-2}\|_q
 = C_q \|f\|_{p,\omega_1}.
$$
In both cases, we see that $\|U(f)\|\le C\|f\|_{\mathcal A}$, so
$U$ is continuous. Clearly $U$ is a homomorphism.

\noi It remains now to show that $U$ is topologically irreducible. This will follow from the fact that closed invariant subspaces of $U$ are invariant under $T$.

\noi Suppose that $Z\subset\ell_1$ is a nonzero invariant subspace of $U$. Take $z\in Z$, $z\ne0$. Let $I_\e$ be the indicator function of $[0,\e]$ and let $\xi_\e=\e^{-1}I_\e$. Then
\begin{align*}
U(\xi_\e)-\Id &= {1\over\e}\int_0^\e e^{xT} dx - {1\over \e} \int_0^\e \Id dx
 = {1\over\e}\int_0^\e (\sum_{n=0}^\infty {(xT)^n\over n!}-\Id) dx
\\&  = {1\over\e}\int_0^\e \sum_{n=1}^\infty {(xT)^n\over n!} dx
  = {1\over\e}\int_0^\e \Big(xT + \sum_{n=2}^\infty {(xT)^n\over n!}\Big) dx.
\end{align*}
Now, assuming that $0<\e<1$, we have
\begin{align*}
\|U(\xi_\e)-\Id -{\e\over 2} T\| &=
   \|{1\over\e}\int_0^\e x \,dx\, T - {\e\over 2} T + {1\over\e}\int_0^\e \sum_{n=2}^\infty {(xT)^n\over n!}\Big) dx\|
  \le {1\over\e} \int_0^\e \sum_{n=2}^\infty {\|(xT)^n\|\over n!} dx
\\& \le {1\over\e} \int_0^\e x^2 \sum_{n=0}^\infty {x^n\over (n+2)!} dx
\le {1\over\e} \int_0^\e x^2 e^x dx \le {1\over\e} \int_0^\e e\, x^2 dx
 = {e\over\e} \cdot {  \e^3\over3} < \e^2.
\end{align*}
Thus, $\e^{-1}\big(U(\xi_\e)-\Id\big) -{1\over 2} T\to 0$, as $\e\to0$, or
$T=\lim\limits_{\e\to0} 2\e^{-1} \big(U(\xi_\e)-\Id\big)$. If now $Z$ is an invariant subspace for $U$, then its closure $\bar Z$ is invariant for $T$, so $\bar Z$ is trivial, and we are done.

\bigskip
\noi {\bf Remarks}. It is clear that this example can be extended to $\R^n$: replace $e^{xT}$ by $e^{\eta(x)T}$, where $\eta$ is a nonzero linear form on $\R^n$. One can show also that $\omega(x)>C\exp(x/\ln x)$, and so the algebras which we construct are not regular.

\bigskip

 \noi Julia Kuznetsova, {\it Unit\'e de Recherche en
Math\'ematiques, Universit\'e du Luxembourg, 6, rue
Coudenhove-Kalergi, L-1359 Luxembourg, Luxembourg,
julia.kuznetsova@uni.lu }

\medskip
\noi Carine Molitor-Braun, {\it  Unit\'e de Recherche en
Math\'ematiques, Universit\'e du Luxembourg, 6, rue
Coudenhove-Kalergi, L-1359 Luxembourg, Luxembourg,
carine.molitor@uni.lu}


\end{document}